\documentclass[12pt]{amsart} 
\usepackage{lmodern}
\usepackage{amsmath}
\usepackage{amssymb}
\usepackage{latexsym}
\def\Nr{{\mathbb N}}

\def\Rr{{\mathbb R}}

\def\Ac{{\mathcal{A}}}
\def\Bc{{\mathcal{B}}}
\def\Cc{{\mathcal{C}}}

\def\Pc{{\mathcal{P}}}

\def\Sc{{\mathcal{S}}}

\def\Vc{{\mathcal{V}}}
\def\Wc{{\mathcal{W}}}

\def\one{{\rm \bf 1}}

\def\dist{\operatorname{dist}}
\def\({\left(}     
\def\){\right)}    
\def\[{\left[}     
\def\]{\right]}

\def\cadlag{c\`adl\`ag{ }}
\newtheorem{definition}{Definition}
\newtheorem{theorem}{Theorem}
\newtheorem{lemma}{Lemma}
\newtheorem{proposition}{Proposition}
\newtheorem{remark}{Remark}
\newtheorem{corollary}{Corollary}

\begin{document}
\author{Freddy Delbaen}
\address{Departement f\"ur Mathematik, ETH Z\"urich, R\"{a}mistrasse
101, 8092 Z\"{u}rich, Switzerland, also Institut f\"ur Mathematik,
Universit\"at Z\"urich, Winterthurerstrasse 190,
8057 Z\"urich, Switzerland}
\email{delbaen@math.ethz.ch}
\title{Monetary Utility Functions on $C_b(X)$ Spaces}
\maketitle

\abstract  We will characterise robust monetary utility functions defined on the space of real valued (bounded) continuous  functions on a Polish space.   \endabstract

\section{Notation and Acknowledgement}
We first recall some definitions from general topology.  A Polish space $X$, is a topological space that is separable and where the topology is metrisable with a complete metric. In this paper we will always use $X$ for a fixed Polish space. No extra assumptions will be put on $X$.  There are no local compactness assumptions as we believe that the most interesting applications are when $X$ is for instance an infinite dimensional Banach space.  The space of real valued bounded continuous functions defined on $X$ is denoted by $C_b(X)$.
\begin{definition} A  function $u\colon C_b(X)\rightarrow \Rr$ is called a concave monetary utility function if it satisfies the following properties
\begin{enumerate}
\item $u$ is a concave function
\item if $f\in C_b(X)$ is nonnegative then $u(f)\ge 0$
\item for $a\in \Rr$ we have $u(f+a)=u(f)+a$ for all $f\in C_b(X)$
\item we say that $u$ is coherent if moreover $u$ is positively homogeneous
\item we say that $u$ satisfies the Fatou property if for every sequence $(f_n)_n$  in $C_b(X)$ so  that  $f_n\downarrow f\in C_b(X)$ pointwise,  $u(f_n)\downarrow u(f)$
\end{enumerate}
A concave monetary utility function  is completely characterised by the acceptability set $\Ac=\{f\mid u(f)\ge 0\}$. For $f\in C_b(X)$ we have $u(f)=\max\{a\in\Rr\mid f-a\in\Ac\}$. 
\end{definition}
\begin{remark}{\rm  The usual definition of the Fatou property is different.  Usually one would require that $u(f)\ge \limsup u(f_n)$ for every uniformly bounded sequence $(f_n)_n$ pointwise converging to $f$.   This is equivalent to the statement that for a uniformly bounded sequence $(f_n)_n$ in $\Ac$ and pointwise converging to $f\in C_b(X)$, one also has $f\in\Ac$.  We only use downward converging sequences and hence we get a more general concept.  But once the main result is proved, one  easily shows, using the Fatou lemma from measure theory, that $u(f)\ge \limsup u(f_n)$ holds for uniformly bounded pointwise converging sequences.
}\end{remark}
From now on we will suppose -- if not otherwise mentioned --  that $u$ is a coherent utility function with the Fatou property. The following proposition is easily proved. We skip the details.
\begin{proposition} $\Ac$ satisfies the following properties:
\begin{enumerate} 
\item $\Ac$ is a convex cone containing the cone of nonnegative functions, $C^+_b(X)$.  In other words $\Ac+C^+_b(X)=\Ac$. If $f\in\Ac$ and $f\le g\in C_b(X)$ then also $g\in\Ac$. 
\item $|u(f)| \le \Vert f\Vert=\sup\{ |f(x)| \mid x\in X\}$. 
\item  If $f_n$ is a sequence in $\Ac$ and $f_n\downarrow f\in C_b(X)$ pointwise, then also $f\in\Ac$.  
\item $\Ac$ is closed for the topology of uniform convergence on $X$.
\item The constant function $-1$ is not in $\Ac$.
\end{enumerate}
\end{proposition}
The main result of the paper is the following
\begin{theorem}  With the notation and properties introduced above there is a weak$^*$ closed convex set $\Sc_0$, of sigma additive probabilities, defined on the Borel sigma algebra of $X$, such that :
$$
\Ac=\{f\mid \mu(f)\ge 0 \text{ for all }\mu \in \Sc_0\}\quad\quad\quad
u(f)=\inf\{\mu(f)\mid \mu\in\Sc_0\}.
$$
\end{theorem}
The difficulty of the proof lies in the fact that in general the dual space of $(C_b(X), \Vert . \Vert)$ is not the space of sigma additive measures on $X$.  For details on the theory of measures on Polish spaces and properties of the weak$^*$ topology (sometimes also called weak topology) the reader can consult \cite{Par}.

We will show how to extend the theorem so that it contains a result for concave monetary utility functions.  The case of coherent utility functions is notationally simpler, so we preferred to prove this case first and then extend it.

When $u$ is defined on the space of all continuous functions on $X$ the Fatou property will be a consequence of the concavity and the monotonicity.  In that case we have

\begin{theorem} Suppose $u\colon C(X)\rightarrow \Rr$ satisfies
\begin{enumerate}
\item $u$ is a concave function that is positively homogeneous
\item $u$ is monetary:  for each $f\in C(X)$ and each $a\in\Rr$ we have $u(f+a)=u(f)+a$
\item $u$ is monotonic: $f\le g$ implies $u(f)\le u(g)$.
\end{enumerate}
Then there is a compact set $L\subset X$ and a convex set, $\Sc$, of probability measures having support contained in $L$, i.e. $\mu\in\Sc$ implies $\mu(L)=1$, such that for all $f\in C(X)$, $u(f)=\inf\{\mu(f)\mid \mu\in\Sc\}$.
\end{theorem}

During the 2014--visit to Shandong University in Jinan, Professor Xinpeng Li asked the author whether he could prove the main result of this paper.  We looked at it but could not find an immediate breakthrough. The problem remained as a cold case on our desk.  We thank Xinpeng Li for the discussions we had on this question.

As possible applications we think of $X$ being a space of trajectories of prices or portfolios. Then a continuous function could be the value of a financial instrument, where the value depends on the trajectory.  Classical examples are when $X$ is the space of continuous functions on a compact interval, endowed with the sup norm. Another example could be the space of \cadlag trajectories endowed with the Skorohod topology. But of course these are just examples of possible applications. The development of such applications is beyond the limited scope of this paper.

The paper is organised as follows.  In the next section we recall some well known results from general topology.  These topological results will then be used in the third section to prove the main theorem. We also give a more constructive proof using some approximation result that allows continuous functions to be seen as limits of uniformly continuous functions.  In section 5 we extend it to the case of concave functions. We show  how to extend monetary utility functions to the space of bounded Borel functions on $X$ in section 6. The extension is in most cases not unique.  Finally we show that for concave monetary functions defined on the set of all continuous functions the algebraic properties imply the Fatou property.

We emphasize that we only use continuity with respect to downward converging sequences.  The case of continuity with respect to increasing pointwise converging sequences is less general since continuity for increasing sequences implies continuity for decreasing sequences. That case was already studied and analysed in the literature, see for instance the paper by Cheridito, Kupper and Tangpi, \cite{CKT}, where general theorems are proved and where a good overview is given. Roughly speaking the stronger convergence implies some kind of compactness.  The precise meaning of this and its relation to some other concepts from functional analysis is outside the scope of this paper, see \cite{convex}.

\section{Some topological results}

In our setting $X$ will be a Polish space.  This means that it is a separable metrisable space for which there exist a  metric for which it becomes complete, i.e. every Cauchy sequence converges.    It is a well known theorem from general topology that says that $X$ can  be embedded as a dense $G_\delta$ set in a compact metric space.  Also the converse holds, if the Polish space $X$ is a dense subspace of a compact topological space $Y$, then $X$ is  a $G_\delta$ set in $Y$.  A $G_\delta$ set is a set that is the countable intersection of open sets. The complement of a $G_\delta$ set is the countable union of closed sets and such a set is called an $F_\sigma$ set. These results were proved by Alexandroff, Hausdorff and Sierpinsky, see the book by Kelley, \cite{Kelley}, p. 208, Problem K. A more general treatment was developed by \v Cech, \cite{Eng}, p. 196. This book also contains a history of this result.\footnote{We thank Professor Eva Colebunders, Free University of Brussels, for bringing Engelking's book to our attention.} The space $C_b(X)$  can be seen as the space of real valued functions on the Stone-\v Cech compactification $\beta X$ of $X$, see \cite{Eng} or \cite{Kelley} for information and properties of $\beta X$. In other words $C_b(X)=C_b(\beta X)=C(\beta X)$. The Stone-\v Cech compactification has the following universal property.  If $f\colon X\rightarrow L$ is a continuous function where $L$ is a compact space, then there is a continuous  extension $\tilde f$ of $f$,  $\tilde f\colon \beta X\rightarrow L$.  There are of course other compactifications of $X$. In particular we can embed $X$ as a dense space of a compact metric space $K$ (in the sequel $K$ will always denote a metrisable compactification of $X$).  Again $X$ is then a $G_\delta$ set in $K$.  So we get the following situation. We have the identity mapping $i\colon X\rightarrow X$. This mapping extends to a continuous function $q\colon \beta X\rightarrow K$.  The following proposition is then an easy consequence.
\begin{proposition} The mapping $q$ (as defined in the previous paragraph) maps $\beta X\setminus X$ onto $K\setminus X$. The set $\beta X\setminus X$ is the  union of a sequence of compact $G_\delta$ sets, $(M_n)_n$.
\end{proposition}
{\bf Proof } The first line of the proposition can be proved as follows.  Suppose $y\in\beta X\setminus X$ is such that $q(y)=x\in X$.  Then of course $q(y)=i(x)=q(x)$. There are  neighbourhoods $U$ of $y$ and $V$ of $x$ in $\beta X$ that are open and disjoint. Because $X$ is dense in $\beta X$, there is a generalised sequence $y_\alpha\in U\cap X$ that converges to $y$. We note that for each $\alpha$, $y_\alpha\notin V$.  Because $q$ is continuous $q(y_\alpha)$ tends to $q(y)=x$.  So there is $\alpha_0$ so that for each $\alpha$ following $\alpha_0$ in the generalised sequence (or net) $y_\alpha=q(y_\alpha)\in V$.  This is a contradiction to $U\cap V=\emptyset$.  Hence $q(\beta X\setminus X)\subset K\setminus X$.  But the compactness of $q(\beta X)$ implies that $q(\beta X)= K$ and hence $q(\beta X\setminus X)=K\setminus X$. 

The set $K\setminus X$ is the union of a sequence of compact sets $(L_n)_n$. These sets are necessarily $G_\delta$ sets since $K$ is a metric space. Since $q$ is continuous the sets $M_n=q^{-1}(L_n)$ are also $G_\delta$ sets.  The first result of the proposition shows that $M_n\subset \beta X\setminus X$ and  $\beta X\setminus X=\cup_{n\ge 1}M_n$.\qed
\begin{corollary} For each set $M_n$ there is an element $\phi_n\in C_b(X)$ such that $0\le \phi_n\le 1$ on $\beta X$ and $M_n=\phi_n^{-1}(1)$.  In particular $\phi_n(x)<1$ for each $x\in X$.  The sets $M_n$ are elements of the sigma-algebra generated by the real valued continuous functions, the so called Baire sigma-algebra $\Bc_0$.  (This is a sigma algebra that is smaller than the Borel sigma-algebra, $\Bc$, which is generated by all closed subsets of $\beta X$).
\end{corollary}\label{phi-n}
{\bf Proof } We only have to give a construction of the functions $\phi_n$. On $K$ we use the distance function to $L_n$ and use the function $k\rightarrow\max(1- \dist(k,L_n),0)$. Then we compose this function with the already defined function $q$.\qed

\begin{remark}{\rm That the sets $M_n$ can be chosen as $G_\delta$ and hence are peak sets is not always mentioned in the literature. It plays a fundamental role in the next sections. An immediate consequence is the obvious fact that $\beta X\setminus X\in \Bc_0$.
}\end{remark}
\begin{remark}{\rm The above reasoning can be generalised to yield the result that when $X$ is densely embedded in the compact space $B$ (in other words $B$ is a compactification of $X$) then $B\setminus X$ is the countable union of compact $G_\delta$ sets. Because this is a little bit outside the scope of this paper we only give a sketch. We first take a countable set of elements of $C(B)$ that generates the topology of $X$ and hence separates the points of $X$.  This is possible since $X$ is a separable metric space.  Then we look at the closed algebra generated by this countable set.  This is a separable Banach algebra equal to $C(K)$ where $K$ is a compact metric space. $K$ is then also a compactification of $X$ and because $X$ is Polish, it is a $G_\delta$ in $K$. It is easily seen that the elements of $B$ define a multiplicative linear form on $C(K)$, hence define elements of $K$.  As above $B\setminus X$ is mapped onto $K\setminus X$ and the proof then continues as before.
}\end{remark}
The dual space of $C_b(X)=C(\beta X)$ is the space of ``regular" sigma-additive measures on $\Bc_0$.  We do not need/use that each such measure can be extended to a regular measure on the Borel sigma algebra. The set of regular sigma additive  nonnegative normalised measures, i.e. probability measures, is denoted by $\Pc$ (or when details are needed $\Pc(X),\Pc(K),\Pc(\beta X),\ldots$) and is equipped with the weak$^*$ topology. With this topology $\Pc(K),\Pc(\beta X)$ are compact spaces.
\section{Coherent monetary utility functions}

We use the notation introduced in the preceding sections. We first note that $\Ac$ is a cone in $C(\beta X)$ that is also closed for the uniform or norm topology.  As in the classical case of the theory of risk measures it is easily seen that there is a weak$^*$ closed convex subset $\Sc\subset \Pc(\beta X)$ such that
\begin{eqnarray*}
\Sc&=&\{\mu\in\Pc(\beta X)\mid \mu(f)\ge 0\text{ for all }f\in\Ac\},\\
\Ac&=&\{f\in C(\beta X)\mid \mu(f)\ge 0\text{ for all }\mu\in \Sc\}.
\end{eqnarray*}
The proof is a simple application of the Hahn-Banach theorem or the bipolar theorem, see \cite{Rudin} for details on functional analysis. The set $\Sc$ is constructed as a base for the polar cone $\Ac^\circ$ of $\Ac$.
The problem we face is that some elements of $\Sc$ might charge the set $\beta X\setminus X$ and hence cannot be seen as sigma additive probabilities on $X$. The proof of the main result reduces to the proof that the subset $\Sc_0\subset \Sc$ of elements $\mu\in\Sc$ such that $\mu(\beta X\setminus X)=0$ is weak$^*$ dense in $\Sc$.
\begin{lemma} For each $1>\epsilon>0$ and each $M_n$ the set
$$
\Vc=\{\mu\in\Sc\mid \mu(M_n)<\epsilon\}
$$
is a dense open set in $\Sc$.
\end{lemma}
{\bf Proof } That the set $\Vc$ is open follows from the definition of the weak$^*$ topology. Indeed $\mu(M_n)<\epsilon$ if and only if there is $k\ge 1$ such that $\mu(\phi_n^k)<\epsilon$. The set $\Vc$ is clearly convex and so is its closure $\Wc$.  Suppose now that $\Wc\neq \Sc$.  Then by the separation theorem \footnote{In case $\Vc=\emptyset$, we can take $f=-1$ and proceed in the same way.} there is an element $f\in C(\beta X)$ such that 
$$
\min \{\mu(f)\mid \mu\in\Sc\} < 0 \le \inf\{\mu(f)\mid \mu\in \Wc\}.
$$
That $0$ is a separating value can be obtained by adding a suitable constant to $f$.  We may suppose that $\Vert f\Vert\le 1$, eventually we normalise $f$. The inequality on the left implies that $f\notin \Ac$. Now take $k$ big enough to ensure that $k\epsilon\ge 1$. For each $m\ge 1$ we define $f_m = f+k\phi_n^m$. Since $\phi_n(x)<1$ for each $x\in X$ we have $f_m\downarrow f$ on $X$. We will now show that $f_m\in \Ac$ for each $m$ and this implies by the hypothesis on $\Ac$ that also $f\in\Ac$, a contradiction. To see that $f_m\in\Ac$ we will calculate $\mu(f_m)$ for each $\mu\in\Sc$.  If $\mu\in\Wc$ then necessarily $\mu(f_m)\ge\mu(f)\ge 0$.  If $\mu\notin\Wc$, then $\mu(M_n)\ge \epsilon$ and $\mu(f_m)\ge -1+k\mu(M_n)\ge -1+k\epsilon\ge 0$. This shows $f_m\in\Ac$ yielding the desired contradiction. \qed

{\bf Proof of the main theorem } The set $\Sc_0$ can be written as
$$
\Sc_0=\bigcap_{n\ge 1, N\ge 1}\left\{\mu\in\Sc\mid \mu(M_n)< \frac{1}{N}\right\}
$$
and is therefore the countable intersection of dense open sets in $\Sc$.  Since $\Sc$ is a compact set, it has the Baire property, meaning that the intersection of a sequence of dense open sets is dense. This shows that
$$
\Ac=\{f\in C(\beta X)=C_b(X)\mid \mu(f)\ge 0\text{ for all }\mu\in\Sc_0\}.
$$
Since $\Sc_0=\{\mu\in\Pc(X)\mid \mu(f)\ge 0\text{ for all }f\in\Ac\}$ it is clear that it is a weak$^*$ closed set in $\Pc(X)$.\qed
\section{Another -- more constructive -- proof of the main theorem}
In this section we will give a more constructive proof of the main theorem.  We will not use the axiom of choice.\footnote{~ That we give such a proof should not be seen as an indication that the author has converted to constructivism}  We will make use of the separation theorem and the compactness of the dual ball of a Banach space in its weak$^*$ topology but since the spaces we will use are separable, the general version of the axiom of choice can be replaced by a countable equivalent. The interested reader can look up such developments in the specialised literature. We do not give details, not even references. The proof we will give also uses some, maybe useful result on the approximation of continuous functions by uniformly continuous functions.

We again use a compactification of $X$ and this time we suppose that the Polish space $X$ is a dense $G_\delta$ of a compact metric space $K$. The metric on $K$ will be denoted by $d$. Bounded continuous functions on $X$ cannot necessarily be extended to continuous functions on $K$.  However we know that a bounded continuous function on $X$ can be extended to a continuous function on $K$ if and only if it is uniformly continuous for $d$.
\begin{lemma} For each element $f\in C_b(X)$ there is a sequence of functions $h_n\in C(K)$ such that $h_n\downarrow f$ on $X$.
\end{lemma}
{\bf Proof } We will construct a sequence of elements $g_n\in C(K)$ such that $g_n\uparrow f$.  By applying the result to $-f$ we then get the statement of the lemma, in fact we then will have proved that there are two sequences in $C(K)$ one approaching $f$ from above and the other approaching $f$ from below.  The proof  uses a convolution technique.  For each $n$ we define
$$
g_n(x)=\inf \{f(y)+n\,d(x,y)\mid y\in X\}.
$$
Clearly $g_n(x)\le f(x)$.  For two points $x_1,x_2$ in $X$ we have
\begin{eqnarray*}
g_n(x_1)&=&\inf\{f(y)+n\,d(x_1,y)\mid y\in X\}\\
&\le& \inf \{f(y)+n\,d(x_2,y) + n\,d(x_1,x_2)\mid y\in X\}\\
&=& g_n(x_2)+n\,d(x_1,x_2).
\end{eqnarray*}
By symmetry we then get $|g_n(x_1)-g_n(x_2)|\le n\, d(x_1,x_2)$, proving that $g_n$ can be extended to an element of $C(K)$. It is also obvious that $g_n\le g_{n+1}$.  We now show that for each $x\in X$, $g_n(x)\uparrow f(x)$.  To estimate $g_n(x)$ we distinguish between $d(x,y)\le \frac{2\Vert f \Vert}{n}$ and $d(x,y)\ge \frac{2\Vert f \Vert}{n}$.  In the latter case we have $f(y)+n\,d(x,y)\ge f(y) + 2\Vert f\Vert\ge \Vert f\Vert\ge f(x)$.  In the former case we have $f(y)+n\,d(x,y)\ge \inf_{d(x,z)\le \frac{2\Vert f \Vert}{n}}f(z)$. We get $g_n(x)\ge \inf_{d(x,y)\le \frac{2\Vert f \Vert}{n}}f(y)$. As $n\rightarrow \infty$ the continuity of $f$ implies that $g_n(x)\rightarrow f(x)$.\qed

Let now $\Ac_u$ be the set of elements of $\Ac$ that are uniformly continuous, identified with their extensions to $K$.  With a slight abuse of notation $\Ac_u=\Ac\cap C(K)$. The set $\Ac_u$ is a cone, closed for the uniform convergence on $K$ and containing the  cone of nonnegative elements of $C(K)$. There exists a convex, weak$^*$ compact set of probabilities on $K$, denoted by $\Sc$, such that
$$
\Ac_u=\{h\in C(K)\mid \mu(h)\ge 0\text{ for all }\mu\in\Sc\}.
$$
The cone $\Ac_u$ satisfies the same sequential closedness property as $\Ac$.  This means that if $h_n$ is a sequence in $\Ac_u$ that converges on $X$ downward to $h\in C(K)$, then also $h\in\Ac_u$. We can now copy the proof of the preceding section and get that $\Sc_0=\{\mu\in\Sc\mid \mu(X)=1\}$ is weak$^*$ dense in $\Sc$.  This gives us
$$
\Ac_u=\{h\in C(K)\mid \mu(h)\ge 0\text{ for all }\mu\in\Sc_0\}.
$$
But the approximation lemma then shows that also
$$
\Ac=\{f\in C_b(X)\mid \mu(f)\ge 0\text{ for all }\mu\in\Sc_0\}.
$$
\qed
\section{Extension to concave monetary utility functions}

For concave monetary utility functions $u$, we will use a  metrisable compactification of $X$, $X\subset K$. The probability measures on $X$ resp. $K$ are denoted by $\Pc(X), \Pc(K)$. As before we put $\Ac_u=\Ac\cap C(K)$.  The Fenchel duality theory then says that $u$ can be calculated as follows.  First we define the conjugate function, which can only take finite values on $\Pc(K)$:
$$
c(\mu)=\sup\{\mu(-f)\mid f\in \Ac_u\}.
$$
The function $c\colon \Pc(K)\rightarrow \overline{\Rr_+}$ is lower semi continuous for the weak$^*$ topology and convex.  This means that the epigraph 
$$\Cc=\{(\mu,t)\mid t\ge c(\mu), t\in\Rr_+, \mu\in\Pc(K)\}$$ is a closed convex set in $\Pc(K)\times \Rr_+$. Observe the trivial fact that $(\mu,t)\in\Cc$ implies that $c(\mu)<+\infty$.
Then the duality result (proved using the Hahn-Banach theorem) says that
\begin{eqnarray*}
u(f)&=& \inf \{\mu(f)+c(\mu)\mid \mu\in\Pc(K)\}=\min \{\mu(f)+c(\mu)\mid \mu\in\Pc(K)\}\\
\Ac_u&=&\{f\in C(K)\mid \mu(f)+c(\mu)\ge 0\text{ for all }\mu\in\Pc(K)\}.
\end{eqnarray*}
\begin{theorem}  Suppose that $\colon C_b(X)\rightarrow \Rr$ satisfies
\begin{enumerate}
\item $u$ is a concave function
\item if $f\in C_b(X)$ is nonnegative then $u(f)\ge 0$
\item for $a\in \Rr$ we have $u(f+a)=u(f)+a$ for all $f\in C_b(X)$
\item $u$ satisfies the Fatou property i.e. for every sequence $(f_n)_n$  in $C_b(X)$ so  that  $f_n\downarrow f\in C_b(X)$ pointwise,  $u(f_n)\downarrow u(f)$,
\end{enumerate}
then the conjugate lower semi continuous convex function $c\colon \Pc(X)\rightarrow \overline{\Rr_+}$ defines $u$:
$$
u(f)=\inf\{\mu(f)+c(\mu)\mid\mu\in\Pc(X)\}
$$
\end{theorem}
{\bf Proof } To prove the theorem it is sufficient to show that the restricted epigraph
$$
\Cc_0=\{(\mu,t)\mid t\ge c(\mu), t\in\Rr_+, \mu\in\Pc(X)\}
$$
is dense in $\Cc$. We proceed as in the previous section using $G_\delta$ sets and the Baire theorem. We first observe that there is a sequence of compact sets $(M_n)_n$ such that $K\setminus X=\cup_n M_n$. We also use the functions $\phi_n$ as defined in Corollary 1. For $\epsilon>0$ we define (again we drop the dependence on $n$ and $\epsilon$):
$$
\Vc=\{(\mu,t)\in\Cc\mid \mu(M_n)<\epsilon\}.
$$
We will show that $\Vc$ is a dense open set in $\Cc$. To do this we look at the closure $\Wc$ in $\Pc(K)\times \Rr$ and suppose that $\Wc\neq\Cc$. We first suppose that $\Vc$ is not empty and will handle this trivial case afterwards. Let us fix an element $(\mu_0,t_0)\in\Cc\setminus \Wc$. Then the separation theorem gives us a couple $(f,\alpha)$ ($f\in C(K), \alpha\in\Rr$) such that
$$
\mu_0(f)+\alpha t_0<0\le \inf\{\mu(f)+\alpha t\mid (\mu,t)\in\Wc\}.
$$
Because $\Wc\neq \emptyset$ we must have that $\alpha\ge 0$, otherwise the right hand side is $-\infty$.  (That $0$ can be used as a separating value can be achieved by adding a constant function to $f$.)  If $\alpha=0$ we can increase $\alpha$ by an arbitrary small positive number so that the inequality remains valid. So we may assume that $\alpha>0$ and we can divide both sides of the equation by $\alpha$. In either case and using $t_0\ge c(\mu_0)$, we get
$$
\mu_0(f)+c(\mu_0)<0\le \inf\{\mu(f)+ c(\mu)\mid (\mu,t)\in\Wc\}.
$$
In case $\Vc$ is empty we take $f=-1$ and proceed as in the case $\alpha=0$ above. 
Elements $(\mu,t)\in\Cc$ that are not in $\Wc$ satisfy $\mu(M_n)\ge\epsilon$.  Now take $k$ big enough to make sure that $-\Vert f\Vert + k\epsilon \ge 0$.  Then the functions $f_N=f+k\phi_n^N$ satisfy $\mu(f_N)+c(\mu)\ge \mu(f_N)\ge 0$ for all $\mu$ with $\mu(M_n)\ge\epsilon$ and $\mu(f_N)+c(\mu)\ge \mu(f)+c(\mu)\ge 0$ if $\mu(M_n)<\epsilon$ because $(\mu,c(\mu))\in\Vc$ in that case. This shows that
 $f_N\in \Ac_u$.  As in section 3 we have $f_N\downarrow f$ on $X$ and hence $f\in \Ac_u$ which is a contradiction. This ends the proof that $\Vc$ is dense in $\Cc$. We must still show that the set $\Vc$ is a $G_\delta$ in $\Cc$.  This is obvious since $\Vc=\Cc\cap\( \{\mu\mid \mu(M_n)<\epsilon\}\times \Rr\)$ and hence is an open set in $\Cc$. The set $\Cc_0$ is the countable intersection of dense open sets of $\Cc$ and since $\Cc$ is a Polish space, it has the Baire property.  Therefore $\Cc_0$ is also a dense set in $\Cc$.\qed
\section{Extension to bounded Borel functions}

A Fatou concave monetary utility function $u\colon C_b(X)\rightarrow \Rr$ can be written using probability measures on the Polish space $X$. As proved above there is a weak$^*$ lower semi continuous function $c:\Pc(X)\rightarrow \overline{\Rr_+}$ such that for $f\in C_b(X)$ we have $u(f)=\inf\{\mu(f)+c(\mu)\mid \mu\in\Pc(X)\}$. This formula can be used for bounded Borel functions and hence we get an extension $\tilde{u}$ of $u$. If the space of bounded Borel functions is denoted by $\Bc_b(X)$  we get
\begin{enumerate}
\item $\tilde{u}\colon \Bc_b(X)\rightarrow \Rr\quad \tilde{u}(f)=\inf\{\mu(f)+c(\mu)\mid \mu\in\Pc(X)\}$
\item $\tilde{u}$ is concave, monetary and $f\le g$ implies $\tilde{u}(f)\le\tilde{u}(g)$
\item it has the Fatou property: if $f_n$ is a sequence in $\Bc_b(X)$ tending downward pointwise to $f\in\Bc_b(X)$ then $\tilde{u}(f_n)\downarrow\tilde{u}(f)$.
\end{enumerate}
However there is no guarantee that the extension is uniquely defined, not even among the Fatou, monetary concave utility functions. To see this we will define two Fatou, coherent utility functions $u_1,u_2$ which agree on $C_b(X)$ but are different on $\Bc_b(X)$. This example even works on perfect compact metrisable spaces.  Let us recall that a set is of first category if it is contained in the countable union of closed sets having empty interior. In other words the complement contains a dense $G_\delta$. We now define the two functions as follows:
\begin{eqnarray*}
 u_1(f)&=&\inf \{f(x)\mid x\in X\}\\
 u_2(f)&=&\sup\{a\mid \{x\mid f(x)\le a\}\text{ is of first category }\}.
\end{eqnarray*}
For the interval $[0,1]$ let  $f$ be the indicator function $\one_J$, of the irrational numbers $J$.  Clearly $u_1(f)=0$ whereas $u_2(f)=1$.  Both functions satisfy the Fatou property, see \cite{general}, example 4.8. For continuous functions we have that both functions agree. The defining set for $u_1$ is the set of all probability measures, whereas for $u_2$ we can only use finitely additive probability measures.
\section{The case $C(X)$}

In this section we assume that $X$ is a non compact Polish space. The space $C(X)$ of all continuous real valued functions defined on $X$,  is usually endowed with the topology of uniform convergence on compact sets.  The dual space is then the space of sigma-additive measures on $X$, having compact support. This suggest that monetary concave utility functions on $C(X)$ are representable by a set of probabilities concentrated on the same compact set. Since the space $C(X)$ is in general not metrisable, we might think that sequential continuity as in the Fatou property, is not sufficient. But if the utility function is defined on such a big space and has some algebraic properties, more is true.
\begin{theorem}  Suppose $u\colon C(X)\rightarrow \Rr$ satisfies
\begin{enumerate}
\item $u$ is a concave function that is positively homogeneous
\item $u$ is monetary:  for each $f\in C(X)$ and each $a\in\Rr$ we have $u(f+a)=u(f)+a$
\item $u$ is monotone: $f\le g$ implies $u(f)\le u(g)$.
\end{enumerate}
Then $u$ also satisfies: for decreasing sequences $f_n\downarrow f$ converging pointwise to a continuous function $f$ we have $u(f_n)\downarrow u(f)$.
\end{theorem}
{\bf Proof } We will show that if $f_n\downarrow f$, $f_n,f\in C(X)$, are such that for  some $\delta>0$ and for all $n$, $u(f_n)> u(f)+2 \delta$  then $u$ cannot take finite values on the whole of $C(X)$. The first step is to replace the sequence $f_n$ by a sequence that tends to $f$ in a stationary way.  For this we look at $g_n=\max(f,f_n-\delta)$. Clearly $u(g_n)\ge u(f_n)-\delta\ge u(f)+\delta$. Since $f_n\downarrow f$, we have that also $g_n\downarrow f$ but there is more.  The open sets $G_n=\{x\mid g_n(x)>f(x)\}$ are decreasing. If $x\in G_n$ then $g_n(x)>f(x)$ but this implies that $f_n(x)>f(x)+\delta$. For the closure this implies that on $\overline{G_n}$, $f_n(x)\ge f(x)+\delta$. Hence if $x\in\cap_n \overline{G_n}$ we would have $\lim_n f_n(x)\ge f(x)+\delta$, a contradiction.  This implies $\cap_n\overline{G_n}=\emptyset$. For every $x$ this means that it can only be an element of $G_n$ for a finite number of natural numbers.  So for every $x$ there is $n_0$ (depending on $x$) such that for all $n>n_0$, we have $g_n(x)=f(x)$, i.e. $g_n(x)\downarrow f(x)$ in a stationary way. For each $n$ we will use the function $(n+1)f-n g_n=f-n(g_n-f)\le f$. From the previous line we deduce that this sequence tends to $f$ in a stationary way. Using that
$$
f=\frac{1}{n+1}\((n+1)f-n g_n\)+\frac{n}{n+1}g_n
$$
and applying concavity we get
$$
u(g_n)-\delta\ge u(f)\ge \frac{1}{n+1} u\( (n+1)f-n g_n\)+\frac{n}{n+1}u(g_n).
$$
This implies
$$
\frac{1}{n+1} u(g_n)-\delta \ge \frac{1}{n+1} u\( (n+1)f-n g_n\)
$$
yielding $u\( (n+1)f-n g_n\)\le -(n+1)\delta +u(g_n)\le -(n+1)\delta +u(g_1)$. We now use the properties of the sequence $G_n$ to prove that
$$
F=\inf\{(n+1)f-n g_n\mid n\ge 1\}
$$
is a continuous function.  The  stationary convergence implies that $F\colon X\rightarrow \Rr$. Let now $x_m\rightarrow x$ be a convergent sequence in $X$.
We look at the compact set $C=\{x_m\mid m\ge 1\}\cup\{x\}$ and observe that
$C\cap \cap_n\overline{G_n}=\emptyset$. The compactness then says that there is $n_0$ with $C\cap \overline{G_{n_0}}=\emptyset$.
For the function $F$ this means:
\begin{eqnarray*}
F(x_m)&=&\inf\{((n+1)f-n g_n)(x_m)\mid n\le n_0-1\}\\
F(x)&=&\inf\{((n+1)f-n g_n)(x)\mid n\le n_0-1\}.
\end{eqnarray*}
Hence $F(x_m)\rightarrow F(x)$. The definition of $F$ and the monotonicity of $u$ give that for each $n$, $u(F)\le -(n+1)\delta+u(g_1)$, contradicting that $u(F)>-\infty$.
\qed
\begin{theorem} Suppose $u\colon C(X)\rightarrow \Rr$ satisfies
\begin{enumerate}
\item $u$ is a concave function that is positively homogeneous
\item $u$ is monetary:  for each $f\in C(X)$ and each $a\in\Rr$ we have $u(f+a)=u(f)+a$
\item $u$ is monotone: $f\le g$ implies $u(f)\le u(g)$.
\end{enumerate}
Then there is a compact set $L\subset X$ and a convex set, $\Sc$, of probability measures having support contained in $L$, i.e. $\mu\in\Sc$ implies $\mu(L)=1$, such that for all $f\in C(X)$, $u(f)=\inf\{\mu(f)\mid \mu\in\Sc\}$.
\end{theorem}
{\bf Proof } The previous theorem shows that $u$ has the Fatou property and hence its restriction to $C_b(X)$ can be represented by a convex set $\Sc$ of probability measures on $X$. So for each $f\in C_b(X)$ we have $u(f)=\inf\{\mu(f)\mid \mu\in\Sc\}$.  But the downward sequential continuity then shows that for $f\in C(X)$, $f$ bounded above, $u(f)=\inf\{\mu(f)\mid \mu\in\Sc\}$. We will now show that the elements of $\Sc$ are concentrated on a compact set. For each $\mu\in\Sc$ we denote by $F_\mu$ the support of $\mu$, this is the smallest closed set that has measure equal to one. We now show that $\cup_{\mu\in\Sc}F_\mu$ is relatively compact. We argue by contradiction and suppose that there is a sequence $\mu_n\in\Sc$, $x_n\in F_{\mu_n}$ so that $(x_n)_n$ has no convergent subsequence and are all different. The proof of the following lemma will ge given later.
\begin{lemma} There is a sequence of functions $\psi_n$ with $0\le \psi_n\le 1$, $\psi_n(x_n)=1$ and such that for each sequence of real numbers $(a_n)_n$ the sum  $\sum_n a_n\psi_n$ converges pointwise to a continuous function.
\end{lemma}
Take now $a_n=\frac{n}{u(-\psi_n)}$.  This is well defined since $u(-\psi_n)<0$, $x_n$ being in the support of $\mu_n$. The function $g=\sum_n a_n\psi_n$ is in $C(X)$ and $u(g)\le -a_n u(-\psi_n)=-n$ by monotonicity. Since this must hold for each $n$, we get a contradiction with $u(g)>-\infty$. Let $L=\overline{\cup_{\mu\in\Sc}F_\mu}$ which is a compact set. For each $f\in C(X)$, bounded above, we now have
$$
u(f)=\inf\{\mu(f)\mid \mu\in\Sc\}=\inf\{\mu(f\one_L)\mid \mu\in\Sc\}.
$$
For $f$ not necessarily bounded above we again use concavity to get the same expression.  Let $\sup\{f(x)\mid x\in L\}=n_0$ and observe that for each $n\ge n_0$, $u(f\wedge n)=u(f\wedge n_0)=\alpha$. Suppose that $u(f)\ge \alpha+\delta$, where $\delta>0$. The function $f-2(f-f\wedge n)=2 (f\wedge n)-f$ is bounded above and on $L$ it is equal to $f$.  Hence $u(2 (f\wedge n)-f)=\alpha$.  Now use concavity of $u$ and $f\wedge n={1\over 2} f +{1\over 2}(2 (f\wedge n)-f)$ to get
$$\alpha=u(f\wedge n)\ge {1\over 2} u(f)+{1\over 2}u(2 (f\wedge n)-f)\ge {1\over 2} (\alpha+\delta)+{1\over 2}\alpha=\alpha +{1\over 2}\delta,$$
a contradiction. This shows that also $u(f)=\alpha$ and hence
$$
u(f)=\inf\{\mu(f)\mid \mu\in\Sc\}=\inf\{\mu(f\one_L)\mid \mu\in\Sc\},
$$
for all $f\ C(X)$.\qed

{\bf Proof of the lemma} Let $d$ be a metric that defines the topology of $X$. We may suppose that $d\le 1$.  For each $J\subset \Nr$, the set $\{x_n\mid n\in J\}$ is closed. Obviously $x_1\notin \{x_n\mid n\ge 2\}$, which is a closed set.  Let $0<\eta_1< {1\over 4} \dist(x_1,\{x_n\mid n\ge 2\})$. Inductively we now define $\eta_k, k\ge 2$ so that  $$0<\eta_k\le \min\(\frac{\eta_{k-1}}{2},{1\over 4}\dist(x_k,\{x_n\mid n\ge {k+1}\}\).$$
For each $n$ let $O_n$ be the open ball around $x_n$ with radius $\eta_n$. It is easy to see that for $n\neq m$, $\overline{O_n}\cap\overline{O_m}=\emptyset$. The sequence of closed sets $\overline{O_n};n\ge 1$ is locally finite, i.e. for each $y\in X$ there is $\epsilon>0$ so that $B(y,\epsilon)$ (the open ball of radius $\epsilon$ around $y$) only intersects a finite number of sets of the form $\overline{O_n}$. Indeed if this would not be true, then for some $y\in X$ and all $\epsilon>0$, the set $\{n\mid B(y,\epsilon)\cap \overline{O_n}\neq \emptyset\}$ would be infinite.  This allows to construct a sequence of points $z_k,k\ge 1$ such that $d(y,z_k)\le 2^{-k}$ and $d(z_k,x_{n_k})\le \eta_{n_k}$ where $n_k\rightarrow \infty$.  Consequently $d(y,x_{n_k})\le 2^{-k}+\eta_{n_k}\rightarrow 0$ implying that $x_{n_k}\rightarrow y$, a contradiction to the choice of the sequence $x_n,n\ge 1$. We know -- and this is easily proved --  that for a locally finite family of closed sets the union is also closed. For each $n$, the function $\psi_n$ can be taken as  $$\psi_n(x)=\frac{\dist(x,O_n^c)\wedge \dist(x_n,O_n^c)}{\dist(x_n,O_n^c)}.$$
For a sequence of real numbers $(a_n)_n$ the sum $h=\sum_n a_n\psi_n$ converges pointwise since the supports of the functions $\psi_n$ are contained in the pairwise disjoint sets $\overline{O_n}$. To see that  $h$ is a continuous function we take a convergent sequence $y_m\rightarrow y$ in $X$ and we will show that $h(y_m)\rightarrow h(y)$. If the point $y$ is not in the closed set $\cup_n\overline{O_n}$, then for $m$ big enough $y_m\notin \cup_n\overline{O_n}$ and hence $h(y_m)=0=h(y)$.  If $y\in \cup_n\overline{O_n}$ then let $k$ be the uniquely defined element with $y\in \overline{O_k}$. But in this case $y\notin \cup_{n\neq k}\overline{O_n}$ and since this is a closed set we again have for $m$ big enough $y_m\notin \cup_{n\neq k}\overline{O_n}$. This means that for all $m$ big enough $h(y_m)=a_k\psi_k(y_m)$ and this converges to $a_k\psi_k(y)=h(y)$.\qed


\begin{thebibliography}{100}

\bibitem{CKT} Cheridito, P., Kupper, M. and Tangpi,L.: Representation of increasing convex functionals with countably additive measures, {\em Studia Math. {\bf 260}, 2021, p. 121--140}

\bibitem{general} Delbaen, F.: Coherent Risk Measures on General Probability Spaces, {\em Advances in Finance and Stochastics, Festschrift in Honour of Dieter Sondermann, Springer, Heidelberg, 2001.}

\bibitem{convex} Delbaen, F.: Convex Increasing Functionals on $C_b(X)$ Spaces, {\em forthcoming}

\bibitem{Eng} Engelking, R.: General Topology -- Revised and completed edition, {\em Sigma Series in Pure Mathematics, Volume 6, Heldermann Verlag, Berlin, 1989.}

\bibitem{Kelley} Kelley, J.L.: General Topology, {\em Graduate Text in Mathematics, Springer, Berlin, 1975, previously Van Nostrand Company, Princeton, 1955.}

\bibitem{Par} Parthasarathy, K.R.: Probability Measures on Metric Spaces, {\em Academic Press, New York, 1967.}


\bibitem{Rudin} Rudin, W.: Functional Analysis, {\em Mcgraw Hill, New York, 1973.}

\end{thebibliography}
\end{document}